\documentclass[11pt]{amsart}
\usepackage{amsthm,latexsym,amssymb,amsmath,mathrsfs,paralist}
\usepackage[all]{xy}

\newtheorem{theorem}{Theorem}[section]

\newcommand{\Log}{\mathrm{Log}}
\newcommand{\Lie}{\mathrm{Lie}}
\newcommand{\id}{\mathrm{id}}

\newcommand{\End}{\mathrm{End}}
\newcommand{\ep}{\varphi_p}

\begin{document}

\title{Some applications of the $p$-adic analytic subgroup theorem}
\author[C. Fuchs and D. H. Pham]{Clemens Fuchs and Duc Hiep Pham}

\begin{abstract}
We use a $p$-adic analogue of the analytic subgroup theorem of W\"ustholz to deduce the transcendence and linear independence of some new classes of $p$-adic numbers. In particular we give $p$-adic analogues of results of W\"ustholz contained in \cite{w3} and generalizations of results obtained by Bertrand in \cite{db5} and \cite{db8}.

\end{abstract}
\maketitle


\section{Introduction}

Transcendence theory is known as one of the fundamental questions in number theory. One of the first outstanding results in transcendence theory is Hermite's proof of the transcendence of $e$ obtained in 1873. His work was extended by Lindemann afterwards in order to prove that $e^\alpha$ is transcendental for any non-zero algebraic number $\alpha$. In particular, this showed that $\pi$ is transcendental, and thus answered the problem of squaring the circle. There are many results further pushing the theory forward. Such important results are due to Gelfond, Schneider, Baker, Masser, Coates, Lang, Chudnovsky, Nesterenko, Waldschmidt, W\"ustholz among others.

Especially, in the 1980's W\"ustholz formulated and proved a very deep and far-reaching theorem, the so-called analytic subgroup theorem (see \cite{w1,w2} or \cite{aw}). It is known as one of the most striking theorems in the complex transcendental number theory with many applications. W\"ustholz in 1983 used his analytic subgroup theorem to give a vast generalization of earlier results on linear independence in the complex domain (see \cite{w3}). To present the result, let $k$ be a non-negative integer, $n_0,n_1,\ldots,n_k$ positive integers and $E_1,\ldots,E_k$ elliptic curves defined over $\overline{\mathbb Q}$.
For each $i$ with $1\leq i\leq k$, denote by $L_i$ the field of endomorphisms of $E_i$, $\wp_i$ the Weierstrass elliptic function associated with $E_i$ and $\Lambda_i$ the lattice of periods of $\wp_i$.
Let $u_1,\ldots,u_{n_0}$ be non-zero elements in $\overline{\mathbb Q}$ and let $u_{i,1},\ldots, u_{i,n_i}$ be elements in the set
$$\mathcal E_i:=\{u\in \mathbb C\setminus \Lambda_i; \wp_i(u)\in\overline{\mathbb Q}\},\quad i=1,\ldots,k.$$
Finally, let $V_0$ be the vector space generated by the logarithms $\log u_1,\ldots, \log u_{n_0}$ over $\mathbb Q$ and $V_i$ the vector space generated by $u_{i,1},\ldots, u_{i,n_i}$ over $L_i$ for $i=1,\ldots,k$.
\emph{If the elliptic curves $E_1,\ldots,E_k$ are pairwise non-isogeneous, then the dimension of the vector space $V$ generated by $1, V_0,\ldots, V_k$ over $\overline {\mathbb Q}$ is determined by
\[\dim_{\overline{\mathbb Q}}V=1+\dim_{\mathbb Q}V_0+\dim_{L_1}V_1+\cdots+\dim_{L_k}V_k.\]}
It can be seen that Baker's theorem on linear forms of logarithms of algebraic numbers is the corresponding statement of the theorem without the elliptic curves $E_1,\ldots, E_k$.

Many of the results above have been transferred into the world of $p$-adic numbers as well. The theory of $p$-adic transcendence was studied by many authors following the line of arguments in the classical complex case. This development was initiated by Mahler in the 1930's (see \cite {ma1,ma2}). He obtained the $p$-adic analogues of the Hermite, Lindemann and Gelfond-Schneider theorems. Afterwards several $p$-adic analogues of known results from the complex case were proved, and new theories of the $p$-adic transcendence were obtained. These contributions are due to Veldkamp (cf. \cite{ve}), G\"unther (cf. \cite{ag}), Brumer (cf. \cite{br}), Adams (cf. \cite{a}), Flicker (cf. \cite{f1,f2}), Bertrand and others. Also a $p$-adic analogue of the W\"ustholz' analytic subgroup theorem has been worked out (see \cite{matev} and the paper \cite{fuchs-pham-past} by the authors which revisited that result and which is based on the results of \cite{fuchs-pham}; observe that this statement has already been mentioned by Bertrand in \cite{bertrand}). Before we go on, we state the $p$-adic analogue of the analytic subgroup theorem of W\"ustholz, since it will be the main tool for the results below. As usual, we denote by $\mathbb C_p$ the completion of $\overline{\mathbb Q_p}$ with respect to the $p$-adic absolute value. Let $G$ be a commutative algebraic group defined over $\overline{\mathbb  Q}$ and $\Lie(G)$ denote the Lie algebra of $G$. Then the set $G(\mathbb C_p)$ of $\mathbb C_p$-points of $G$ is a Lie group over $\mathbb C_p$ whose Lie algebra is given by $\Lie(G(\mathbb C_p))=\Lie(G)\otimes_{\overline{\mathbb Q}}\mathbb C_p$. It is known that there is the $p$-adic logarithm map \[\log_{G(\mathbb C_p)}: G(\mathbb C_p)_f\rightarrow \Lie(G(\mathbb C_p))\] where $G(\mathbb C_p)_f$ is the set of $x\in G(\mathbb C_p)$ for which there exists a strictly increasing sequence $(n_i)$ of integers such that $x^{n_i}$ tends to the unity element of $G(\mathbb C_p)$ as $i$ tends to $\infty$ (see \cite[Chapter III, 7.6]{Bourbaki}). We denote by $G_f(\overline{\mathbb Q}):=G(\mathbb C_p)_f\cap G(\overline{\mathbb Q})$ the set of algebraic points of $G$ in $G(\mathbb C_p)_f$. The following statement is the $p$-adic analytic subgroup theorem:
\begin{theorem}\label{past}
Let $G$ be a commutative algebraic group of positive dimension defined over $\overline{\mathbb Q}$ and let $V\subseteq \Lie(G)$ be a non-trivial
$\overline{\mathbb Q}$-linear subspace. For any $\gamma\in G_f(\overline{\mathbb Q})$ with
$0\neq\log_{G(\mathbb C_p)}(\gamma)\in V\otimes_{\overline{\mathbb Q}}\mathbb C_p$
there exists an algebraic subgroup $H\subseteq G$ of positive dimension defined over $\overline{\mathbb Q}$ such that {\rm Lie}$(H)\subseteq V$ and $\gamma\in H(\overline{\mathbb Q})$.
\end{theorem}
In this paper we give some applications of the $p$-adic analytic subgroup theorem. It is known that in the complex domain many results on linear independence and transcendence are deduced from W\"ustholz' analytic subgroup theorem. We shall show that it is also possible to obtain $p$-adic analogues of such results by using the $p$-adic analytic subgroup theorem. We will not give proofs of results which already appeared in the literature but instead concentrate on applications which do not exist in the literature so far.

\section{Results}

We start by giving a $p$-adic analogue of the result of W\"ustholz in \cite{w3} mentioned in the introduction. With the notations as above, but in the $p$-adic domain we let $\wp_{p,i}$ be the (Lutz-Weil) $p$-adic elliptic function associated with $E_i$ (it is given by the same formal power series as the Weiterstrass $\wp$-function in the complex case, however the domain on which it is analytic depends on $p$), and $v_{i,1},\ldots, v_{i,n_i}$ elements in the set
\[\mathcal E_{p,i}:=\{0\}\cup\{u\in \mathscr D_{p,i}\setminus\{0\}; \wp_{p,i}(u)\in\overline{\mathbb Q}\}\]
where $\mathscr D_{p,i}$ is the $p$-adic domain of $E_i$ for $i=1,\ldots,k$ (see \cite{we} and \cite{lu}). We have the following result.

\begin{theorem}\label{n1}
Let $V_0$ be the vector space generated by the $p$-adic logarithms $\Log_p(u_1),\ldots,\Log_p(u_{n_0})$ over $\mathbb Q$, and $V_i$ the vector space generated by $v_{i,1},\ldots, v_{i,n_i}$ over $L_i$ for $i=1,\ldots,k$. If the elliptic curves $E_1,\ldots,E_k$ are pairwise non-isogeneous,
then the dimension of the vector space $V$ generated by $1, V_0,\ldots, V_k$ over $\overline {\mathbb Q}$ is determined by
$$\dim_{\overline{\mathbb Q}}V=1+\dim_{\mathbb Q}V_0+\dim_{L_1}V_1+\cdots+\dim_{L_k}V_k.$$
\end{theorem}

Following \cite{w3} we turn to abelian varieties. Let $A$ be an abelian variety defined over $\overline{\mathbb Q}$ and denote by $\frak a$ the Lie algebra of $A$. It is well-known that $A(\mathbb C_p)$ has naturally the structure of a Lie group over $\mathbb C_p$. By \cite [Chapter III, 7.2, Prop. 3]{Bourbaki}, there is an open subgroup $\mathcal A_p$ of $\Lie(A(\mathbb C_p))=\frak a\otimes_{\overline{\mathbb Q}}\mathbb C_p$ on which the $p$-adic exponential map $\exp_{A(\mathbb C_p)}$ is defined and
\[\exp_{A(\mathbb C_p)}: \mathcal A_p\rightarrow \exp_{A(\mathbb C_p)}(\mathcal A_p)\subset A(\mathbb C_p)\]
is an isomorphism. Its inverse is the restriction of the $p$-adic logarithmic map to $\exp_{A(\mathbb C_p)}(\mathcal A_p)$. We denote by $\mathcal A_{\overline {\mathbb Q}}$ the set
\[\{\alpha\in \mathcal  A_p; \exp_{A(\mathbb C_p)}(\alpha)\in A(\overline{\mathbb Q})\}.\]
Let $R$ be a ring operating on $\frak a$. We say that elements $u_1,\ldots,u_m$ in $\frak a$ are linearly independent over $R$ if whenever we have a relation $r_1u_1+\cdots+r_mu_m=0$ with $r_1,\ldots,r_m$ in $R$, then it implies that $r_1=\cdots=r_m=0.$
Finally, let $\End(A)$ denote the ring of endomorphisms of $A$ and $\End(\frak a)$ the ring of endomorphisms of $\frak a$ (defined over $\overline{\mathbb Q}$). We prove the following theorem which again is the $p$-adic analogue of a result of W\"ustholz in \cite{w3}.

\begin{theorem}\label{n2}
Assume that $A$ is a simple abelian variety. Let $m$ be a positive integer, and $u_1,\ldots, u_m$ elements of $\mathcal A_{\overline{\mathbb Q}}$ such that they are linearly independent over $\End(A)$. Then $u_1,\ldots,u_m$ are linearly independent over $\End(\frak a)$.
\end{theorem}

We now generalize a result of Bertrand from $p$-adic transcendence. To describe the results below, as usual, we identify the Lie algebra $\frak a$ with the vector space $\overline{\mathbb Q}^n$, here $n:=\dim A$. Under this identification it follows that $\mathcal A_{\overline{\mathbb Q}}$ is a subset of
\[\Lie(A(\mathbb C_p))=\frak a\otimes_{\overline{\mathbb Q}}\mathbb C_p=\overline{\mathbb Q}^n\otimes_{\overline{\mathbb Q}}\mathbb C_p=\mathbb C_p^n.\]
We say that $A$ is an abelian variety of {RM type} (or with {real multiplication}) if there are a totally real number field $F$ of degree $n$ and an embedding $F\hookrightarrow \End(A)\otimes_{\mathbb Z}\mathbb Q$ (hence an abelian variety of CM type is also of RM type). The following result was given by Bertrand in 1979 (see \cite{db8}): \emph{If $A$ is a simple abelian variety of RM type defined over $\overline{\mathbb Q}$, then all coordinates of $u$ are transcendental for any non-zero element $u$ in $\mathcal A_{\overline{\mathbb Q}}$.} We shall use the $p$-adic analytic subgroup theorem to give a generalization of this theorem for general simple abelian varieties. That is, we prove the following result.

\begin{theorem}\label{n3}
If $A$ is a simple abelian variety of positive dimension defined over $\overline{\mathbb Q}$, then all coordinates of $u$ are transcendental for any non-zero element $u$ in $\mathcal A_{\overline{\mathbb Q}}$.
\end{theorem}

It is known that there are homomorphisms from the additive group $\mathbb C_p$ to the multiplicative group $\mathbb C_p^*$ extending the $p$-adic exponential function (see \cite[Chapter 5, 4.4]{rm}). Let $\ep$ be such a homomorphism. We get the following result which can be seen as a generalization of Corollary 2 in \cite{db5}.

\begin{theorem}\label{n4}
If $A$ is a simple abelian variety defined over $\overline{\mathbb Q}$ of dimension $n>0$, then the elements $\ep(\alpha u_1),\ldots,\ep(\alpha u_n)$ are transcendental for any non-zero element $u=(u_1,\ldots,u_n)$ in $\mathcal A_{\overline{\mathbb Q}}$ and for any non-zero element $\alpha\in\overline{\mathbb Q}$.
\end{theorem}

In the next section we give the proofs of the statements as consequence of the $p$-adic analytic subgroup theorem.

\section{Proofs}

\subsection{Proof of Theorem \ref{n1}}

The proof is very similar to the one given in \cite{w3}. First of all, one may assume without loss of generality that $n_0=\dim_{\mathbb Q}V_0$ and $n_i=\dim_{L_i}V_i$ for $0\leq i\leq k$. Let $r_{i}$ be the $p$-adic valuation of $u_{i}$, and define $v_{i}:=p^{-r_{i}}u_{i}$ for $i=1,\ldots,n_0$. Then $v_{i}$ are algebraic numbers in $\mathbb U(1):=\{x\in\mathbb C_p; |x|_p=1\}$ and
$$\Log_p(v_i)=\Log_p(u_i),\quad \forall i=1,\ldots,n_0.$$
We therefore have to show that the elements
$$1, \Log_p(v_1),\ldots,\Log_p(v_{n_0}), u_{1,1},\ldots, u_{1,n_1},\ldots, u_{k,1},\ldots, u_{k,n_k}$$
are linearly independent over $\overline{\mathbb Q}$. In fact, assume on the contrary that they are not linearly independent over $\overline{\mathbb Q}$. This means that there exists a non-zero linear form $l$ in $n:=1+n_0+\cdots+n_k$ variables with coefficients in $\overline{\mathbb Q}$ such that
$$l(1,\Log_p(v_1),\ldots,\Log_p(v_{n_0}),u_{1,1},\ldots,u_{1,n_1},\ldots,u_{k,1},\ldots,u_{k,n_k})=0.$$
Let $V$ be the $\overline{\mathbb Q}$-vector space defined by
$$V:=\{v\in \overline{\mathbb Q}^{n}; l(v)=0\},$$
and consider the commutative algebraic group $G$ given by
$$G=\mathbb G_a\times \mathbb G_m^{n_0}\times E_1^{n_1}\times\cdots\times E_k^{n_k}.$$
Then $G$ is defined over $\overline{\mathbb Q}$ and the Lie algebra
$$\Lie(G)=\Lie(\mathbb G_a)\times \Lie(\mathbb G_m)^{n_0}\times \Lie(E_1)^{n_1}\times\cdots\times \Lie(E_k)^{n_k}$$
which is identified with $\overline{\mathbb Q}^{1+n_0+\cdots+n_k}=\overline{\mathbb Q}^n$.
This shows that
$$\Lie(G(\mathbb C_p))=\Lie(G)\otimes_{\overline{\mathbb Q}}\mathbb C_p=\mathbb C_p^{n}.$$
We write abbreviately $\exp_i$ and $\log_i$ for the $p$-adic exponential and logarithm map associated with $E_i$ respectively for each $i=1,\ldots,n$. One has
$$G(\mathbb C_p)_f=\mathbb C_p\times \mathbb U(1)^{n_0}\times (E_1(\mathbb C_p))^{n_1}\times (E_k(\mathbb C_p))^{n_k}$$
and the $p$-adic logarithm map $\log_{G(\mathbb C_p)}: G(\mathbb C_p)_f \rightarrow \mathbb C_p^n$
is determined by
$$(\id_{\mathbb C_p}, (\Log_p)^{n_0},\log_1^{n_1},\ldots,\log_k^{n_k}).$$
Let $\gamma\in G(\mathbb C_p)_f$ be the algebraic point given by
$$\big(1, v_1,\ldots, v_{n_0}, \exp_1(u_{1,1}),\ldots,\exp_1(u_{1,n_1}),\ldots,\exp_k(u_{k,1}),\ldots,\exp_k(u_{k,n_k})\big).$$
One gets
$$\log_{G(\mathbb C_p)}(\gamma)=(1,\Log_p(v_1),\ldots, \Log_p(v_{n_0}), u_{1,1},\ldots,u_{1,n_1},\ldots,u_{k,1},\ldots,u_{k,n_k})$$
which is a non-zero element in $V\otimes_{\overline{\mathbb Q}}\mathbb C_p$. We apply the $p$-adic analytic subgroup theorem to $G,V$ and $\gamma$ to obtain an algebraic subgroup $H$ of $G$ of positive dimension defined over $\overline{\mathbb Q}$ such that $\gamma$ is in $H(\overline{\mathbb Q})$ and $\Lie(H)$ is contained in $V$. Let $H^0$ be the connected component of $H$ then
$$H^0=H_{-1}\times H_0\times\cdots\times H_k$$
where $H_{-1}$ is an algebraic subgroup of $\mathbb G_a$, $H_0$ is an algebraic subgroup of $\mathbb G_m$ and $H_i$ is an algebraic subgroup of $E^{n_i}_i$ for $1\leq i\leq k$.
We have
$$\dim H^0=\dim H=\dim_{\overline{\mathbb Q}}\Lie(H)\leq \dim_{\overline{\mathbb Q}}V=n-1.$$
Since $\gamma\in H$ and the first coordinate of $\gamma$ is $1$ it follows that $H_{-1}=\mathbb G_a$. This means that at least one of the algebraic groups $H_0,\ldots,H_k$ is contained in the corresponding factor with positive codimension. If this happens for the algebraic group $H_0$, i.e. $H_0$ is a proper algebraic subgroup of $\mathbb G_m$, then the Lie algebra of $H_0$ is defined by non-zero linear forms with integer coefficients. This means that the elements $\Log_p(v_1),\ldots,\Log_p(v_{n_0})$ are linearly dependent over $\mathbb Q$ since $(v_1,\ldots,v_{n_0})\in H_0$, i.e. $\dim_{\mathbb Q}V_0<n_0$, a contradiction. Hence there must be at least one $i\in\{1,\ldots,k\}$ such that
$$\dim H_i\leq \dim(E_i^{n_i})-1=n_i-1.$$
Using the same arguments as in the proof of \cite[Theorem 6.2]{aw}, we conclude that the elements $u_{i,1},\ldots,u_{i,n_i}$ are linearly dependent over $L_i$. This is, $\dim_{L_i}V_i<n_i$, a contradiction, and the theorem follows.
\hfill$\square$

\subsection{Proof of Theorem \ref{n2}}

Suppose that there are elements $r_1,\ldots,r_m$ not all zero in $\End(\frak a)$ such that
$$r_1u_1+\cdots+r_mu_m=0.$$
Let $V$ be the vector space defined by
$$V:=\{v=(v_1,\ldots,v_m)\in \frak a^n; r_1v_1+\cdots+r_mv_m=0\}.$$
Then $V$ is a $\overline{\mathbb Q}$-linear subspace of $\frak a^m$. We consider the abelian variety $G:=A^m$ and the element $\gamma:=(\exp_{A(\mathbb C_p)}(u_1),\ldots,\exp_{A(\mathbb C_p)}(u_m))$.
Then $\gamma$ is an algebraic point of
$$G(\mathbb C_p)_f=A(\mathbb C_p)^m_f=A(\mathbb C_p)^m,$$
and $\log_{A(\mathbb C_p)}(\gamma)=(u_1,\ldots,u_m)$ is a non-zero element in $V\otimes_{\overline{\mathbb Q}}\mathbb C_p$.
The $p$-adic analytic subgroup theorem then shows that there exists an algebraic subgroup $H$ of $G$ of positive dimension defined over $\overline{\mathbb Q}$ such that $\gamma\in H(\overline{\mathbb Q})$ and $\Lie(H)$ is a subspace of $V$. Since $A$ is simple, $H$ is isogeneous to $A^k$ for a certain integer $k<m$. One can therefore define a projection $\pi$ from $G$ onto $H$. We get the corresponding tangent map $d\pi: \Lie(G)\rightarrow \Lie(H)$. Since $\gamma\in H(\overline Q)$ it follows that the point $(u_1,\ldots,u_m)=\log_{A(\mathbb C_p)}(\gamma)=\log_{H(\mathbb C_p)}(\gamma)$ belongs to $\Lie(H)$.  On the other hand, we may identify $\Lie(G)$ and $\Lie(H)$ with $\frak a^m$ and $\frak a^k$, respectively. Then the point $(u_1,\ldots, u_m)$ is in the kernel of the linear map $\id_{\Lie(G)}-d\pi$ which can be written as an $m\times m$ matrix with entries in $\End(A)$ (since the algebra of endomorphisms $\End(G)$ of $G$ is represented on the $\Lie(G)$ by the matrix algebra $M_m(\End(A))$). In other words, the image of $(u_1,\ldots, u_m)$ under this matrix is zero. Note that this matrix is non-zero (since $k<m$), hence there is at least one  column of its is non-zero. We have thus shown that there is a non-trivial dependence relation over $\End(A)$ among the elements $u_1,\ldots, u_m$, or equivalently, $u_1,\ldots, u_m$ are linearly dependent over $\End(A)$. This contradiction proves the theorem.
\hfill$\square$

\subsection{Proof of Theorem \ref{n3}}
Denote by $n$ the dimension of $A$. Suppose on the contrary that there is a non-zero element $u=(u_1,\ldots,u_n)\in \mathcal A_{\overline{\mathbb Q}}$ with $u_i$ algebraic over $\overline{\mathbb Q}$ for some $i\in\{1,\ldots,n\}$. Let $G=\mathbb G_a\times A$ be the direct product of the additive group $\mathbb G_a$ with $A$. Then $G$ is commutative and defined over $\overline{\mathbb Q}$ with
$$\Lie(G)=\Lie(\mathbb G_a)\times \Lie(A)=\overline{\mathbb Q}^{n+1}$$
This gives
$$\Lie(G(\mathbb C_p))=\Lie(G)\otimes_{\overline{\mathbb Q}}\mathbb C_p=\mathbb C_p^{n+1}.$$
We have $G(\mathbb C_p)_f=\mathbb C_p\times A(\mathbb C_p)$ and the $p$-adic logarithm map $\log_{G(\mathbb C_p)}$ is $\id_{\mathbb C_p}\times \log_{A(\mathbb C_p)}$. Let $V$ be the $\overline{\mathbb Q}$-vector space defined by
$$\{(v_0,\ldots,v_n)\in \overline{\mathbb Q}^{n+1}; v_0-v_i=0\}$$
and $\gamma$ the algebraic point $(u_i,\exp_{A(\mathbb C_p)}(u))$. We observe that $\gamma$ is non-zero since $u$ is non-zero, and furthemore we have
$$\log_{G(\mathbb C_p)}(\gamma)=\Big(\id_{\mathbb C_p}(u_i), \log_{A(\mathbb C_p)}\big(\exp_{A(\mathbb C_p)}(u)\big)\Big)=(u_i,u_1,\ldots,u_n)$$
is an element in $V\otimes_{\overline{\mathbb Q}}\mathbb C_p.$
Applying the $p$-adic analytic subgroup theorem, we obtain an algebraic subgroup $H$ of $G$ of positive dimension defined over $\overline{\mathbb Q}$ such that $\gamma\in H(\overline{\mathbb Q})$ and $\Lie(H)\subseteq V$. Note that $V$ is proper in $\Lie(G)$, and this shows that $H$ is proper in $G$. Since the abelian variety $A$ is simple, $H$ must be either of the form $\mathbb G_a\times \{e\}$ or $\{0\}\times A$ where $e$ is the identity element of $A$. If $H=\mathbb G_a\times \{e\}$ then $\exp_{A(\mathbb C_p)}(u)=e$, i.e. $u=0$, a contradiction. If $H=\{0\}\times A$ then the Lie algebra $\Lie(H)=\{0\}\times \overline{\mathbb Q}^n$. This contradicts with the condition $\Lie(H)\subseteq V$, and the theorem is proved.
\hfill$\square$

\subsection{Proof of Theorem \ref{n4}}
Assume on the contrary that there exists $i\in\{1,\ldots,n\}$ such that $\ep(\alpha u_i)\in\overline{\mathbb Q}$. There is a sufficiently large positive integer $r$ such that $w:=p^r\alpha u_i\in B(r_p)$, where $B(r_p)$ denotes the ball $\{x\in \mathbb C_p; |x|_p<r_p\}$ with $r_p:=p^{-1/(p-1)}$. Hence
$$e_p(w):=\sum_{k\geq 1}\frac{w^k}{k!}=\ep(w)=\ep(p^r\alpha u_i)=\ep(\alpha u_i)^{p^r}$$
is also in $\overline{\mathbb Q}$. Let $G=\mathbb G_m\times A$ be the direct product of the multiplicative group $\mathbb G_m$ with $A$. Then $G$ is commutative and defined over $\overline{\mathbb Q}$ with
$$\Lie(G)=\Lie(\mathbb G_m)\times \Lie(A)=\overline{\mathbb Q}^{n+1}.$$
This implies that
$$\Lie(G(\mathbb C_p))=\Lie(G)\otimes_{\overline{\mathbb Q}}\mathbb C_p=\mathbb C_p^{n+1}.$$
We have
$$G(\mathbb C_p)_f=\mathbb G_m(\mathbb C_p)_f\times A(\mathbb C_p)_f=\mathbb U(1)\times A(\mathbb C_p),$$
and the $p$-adic logarithm map $\log_{G(\mathbb C_p)}$ is the product $\Log_p\times \log_{A(\mathbb C_p)}$. Let $V$ be the $\overline{\mathbb Q}$-vector space defined by
$$\{(v_0,\ldots,v_n)\in \overline{\mathbb Q}^{n+1}; v_0-p^r\alpha v_i=0\}$$
and $\gamma$ the algebraic point $(e_p(w),\exp_{A(\mathbb C_p)}(u))$. We see that $\gamma$ is non-zero since $u$ is non-zero, and furthermore
$$\log_{G(\mathbb C_p)}(\gamma)=\Big(\log_p\big(e_p(w)\big), \log_{A(\mathbb C_p)}\big(\exp_{A(\mathbb C_p)}(u)\big)\Big)=(p^r\alpha u_i,u_1,\ldots,u_n)$$
is a non-zero element in $V\otimes_{\overline{\mathbb Q}}\mathbb C_p.$
Using the $p$-adic analytic subgroup theorem, we obtain an algebraic subgroup $H$ of $G$ of positive dimension defined over $\overline{\mathbb Q}$ such that $\gamma\in H(\overline{\mathbb Q})$ and $\Lie(H)\subseteq V$. Note that $V$ is proper in $\Lie(G)$, and this shows that $H$ is proper in $G$. Since the abelian variety $A$ is simple, $H$ must be either of the form $\mathbb G_m\times \{e\}$ or $\{1\}\times A$ where $e$ is the identity element of $A$. If $H=\mathbb G_m\times \{e\}$ then $\exp_{A(\mathbb C_p)}(u)=e$, i.e. $u=0$, a contradiction. If $H=\{1\}\times A$ then the Lie algebra $\Lie(H)=\{0\}\times \overline{\mathbb Q}^n$. This contradicts with the condition $\Lie(H)\subseteq V$, and the theorem is proved.
\hfill$\square$

\section{Acknowledgments}

The first author was supported by Austrian Science Fund (FWF): P24574. The second author was supported by grant PDFMP2\_122850 funded by the Swiss National Science Foundation (SNSF).

\vspace*{0.5cm}

\noindent{\sc Clemens Fuchs}\\
Department of Mathematics\\ University of Salzburg\\ Hellbrunnerstr. 34\\ 5020 Salzburg\\ Austria\\
Email: {\sf clemens.fuchs@sbg.ac.at}\\

\noindent{\sc Duc Hiep Pham}\\
University of Education\\Vietnam National University, Hanoi\\ 144 Xuan Thuy, Cau Giay, Hanoi\\ Vietnam\\
Email: {\sf phamduchiepk6@gmail.com}

\end{document}